\documentclass[12pt]{amsart}

\usepackage{hyperref}
 \hypersetup{ colorlinks=true, linkcolor=blue, filecolor=magenta, urlcolor=cyan, }

\usepackage{graphics}
\usepackage{graphicx}
\usepackage{epsfig}
\usepackage{amsmath}
\usepackage{graphics}
\usepackage{graphicx}
\usepackage{epstopdf}
\usepackage{multirow}
\usepackage{booktabs}
\usepackage{subfigure}
\usepackage{ifpdf}
\usepackage{color}
\usepackage{amssymb}
\usepackage{cite}
\usepackage{amssymb}
\usepackage{amsthm}
\usepackage{amsmath}
\usepackage{matlab-prettifier}
\newtheorem*{ntn}{Notations}

\newtheorem{dfn}{Definition}[section]
\newtheorem{thm}{Theorem}[section]

\usepackage{lineno}
\numberwithin{equation}{section}
\newcommand{\rank}{\operatorname{rank}}

\begin{document}


\bibliographystyle{plain}

\title{Finding Tensor Decompositions with Sparse Optimization}

\author{Taehyeong Kim\textsuperscript{1}}
\email{th\_kim@pusan.ac.kr}

\author{Jeong-Hoon Ju\textsuperscript{1}}
\email{jjh793012@naver.com}

\author{Yeongrak Kim\textsuperscript{1,2}}
\email{yeongrak.kim@pusan.ac.kr}

\address{\textsuperscript{1}Department of Mathematics, Pusan National University, 
2 Busandaehak-ro 63beon-gil, Geumjeong-gu, 46241 Busan, Republic of Korea 
\newline
\newline \indent
\textsuperscript{2} Institute of Mathematical Science, Pusan National University, 2 Busandaehak-ro 63beon-gil, Geumjeong-gu, 46241 Busan, Korea}

\begin{abstract}
In this paper, we suggest a new method for a given tensor to find CP decompositions using a less number of  rank $1$ tensors. 
The main ingredient is the Least Absolute Shrinkage and Selection Operator (LASSO) by considering the decomposition problem as a sparse optimization problem. 
As applications, we design experiments to find some CP decompositions of the matrix multiplication and determinant tensors. 
In particular, we find a new formula for the $4 \times 4$ determinant tensor as a sum of $12$  rank $1$ tensors.
\end{abstract}



\maketitle

\section{Introduction}

A tensor is a multilinear map from a product of vector spaces to a vector space. Using ordered bases of the vector spaces, one can represent a multilinear map as a multi-dimensional array and vice versa, similar to the case of linear transformations and matrices. Therefore, a tensor, as a generalization of matrices, is an attractive object from both theoretical and practical points of view. Complex data structures, such as images, videos, and social networks can be represented by tensors as well as classical representations using matrices. There are several tools to analyze matrices such as eigenvalue decomposition or singular value decomposition (SVD), however, these tools are not very well established for tensors at this moment. 

Canonical Polyadic (CP) decomposition is one of the methods to analyze the structure of a given tensor. Roughly speaking, a CP decomposition  expresses a tensor into a sum of rank $1$ tensors \cite{hitchcock1927expression, landsberg2012tensors, kolda2009tensor}. From the viewpoint of tensors as multi-dimensional arrays of data, a CP decomposition is helpful to extract prominent features of the data and to take a low-rank approximation of the given data. For instance, CP decomposition has been widely used in chemometrics and signal processing for extracting the underlying components of a mixture, such as the individual chemical compounds in a sample or the different frequency components in a signal \cite{kolda2009tensor, appellof1981strategies, harshman1970foundations}.


However, it is extremely hard to find more efficient CP decompositions, that is, CP decompositions consist of a lesser number of rank $1$ tensors. The smallest number of rank $1$ tensors required in CP decompositions of a given tensor is called the tensor rank, which generalizes the notion of the matrix rank. Therefore, the tensor rank of a given tensor exhibits the complexity of the given tensor, both as a multi-dimensional array and a multilinear map. Note that the number of summands in a CP decomposition of the given tensor provides an upper bound of the tensor rank, and hence searching for an efficient CP decomposition is also crucial to determine the tensor rank.

In this paper, we propose a point of view so that we may attack these kinds of problems via sparse optimization problems. Indeed, a CP decomposition represents a given tensor as a linear combination of rank $1$ tensors, hence, if we can reduce the number of summands among the candidate rank $1$ tensors, then it will give rise to a more efficient CP decomposition. There are two famous machine learning techniques for solving sparse optimization problems, namely, Sparse Identification of Nonlinear Dynamics (SINDy) and Sparse Dictionary Learning (SDL). 
SINDy is used to identify the underlying mathematical equations that govern a dynamical system from time series data \cite{brunton2016discovering}. 
It is particularly useful for identifying sparse and low-dimensional models from high-dimensional and noisy data. 
SDL is used to find a sparse and optimized solution which leads to the best possible representation of the given data as a linear combination of basic elements in the dictionary \cite{kreutz2003dictionary, mairal2010online}. 
Both methods are widely used in various research fields including signal processing \cite{engan2007family}, computer vision \cite{wright2010sparse}, and data analysis \cite{brunton2016discovering}.

Let us make a closer observation of these methods. Both SINDy and SDL require a sparse regression method, and one of the most popular methods is so-called the Least Absolute Shrinkage and Selection Operator (LASSO). 
LASSO is a non-convex optimization method that can perform both the feature selection and the regularization by the $l_1$-penalty term on the coefficients \cite{santosa1986linear, tibshirani1996regression}. 
Thanks to the $l_1$-penalty term, LASSO can be used in SDL and SINDy to capture the set of sparse elements in the dictionary that best represents the input data \cite{donoho2003optimally, donoho2006most}.
To be precise, LASSO aims to find $X$ that minimizes the loss function $L$, defined as
\begin{equation}\label{LASSO}
	L = \|Y - DX\|_{F} + \lambda\|X \|_1 ,
\end{equation}
for the target equation $Y\approx DX$ where $Y$ is the observation data, $D$ is the candidate data, $\lambda$ is the regularization coefficient, and $\|\cdot \|_{F}$ is the Frobenius norm. 
Hence, LASSO traces an optimal solution that has a small error (corresponding to $\|Y-DX\|_F$) and is sparse (corresponding to $\lambda \|X \|_1$) as possible. A solution $X$ becomes sparser as $\lambda$ increases, and more accurate as $\lambda$ decreases. 
In practice, cross-validation is used to find an optimal choice for the parameter $\lambda$, or alternatively, we may select an optimal $\lambda$ by graphing the change of coefficients with respect to the value for $\lambda$ \cite{tibshirani1996regression}. Therefore, if we select a nice dictionary $D$ consisting of rank $1$ tensors, we expect a sparse solution to express the given tensor $Y$ as a linear combination of rank $1$ tensors in $D$ by solving group sparse optimization problem.

The structure of the paper is as follows. 
In Section \ref{Sect:Preliminaries}, we review basic notions about the rank of a tensor, SINDy, and SDL. 
In Section \ref{Sect:Method}, we describe our method for how to obtain a CP decomposition based on the inspiration from SINDy and SDL.
As applications, we apply our methods to various examples looking for rank and CP-decomposition of tensors in Section \ref{Sect:Result}. In detail, we test with a matrix multiplication tensor $M_{\langle 2 \rangle}$ and determinant tensors. In particular, on Theorem \ref{MainThm} we find a new formula of the determinant for $4\times 4$ matrices, which significantly improves an upper bound for the tensor rank of the $4\times 4$ determinant tensor.

\section*{Acknowledgement}
J.-H. J. and Y. K. are supported by the Basic Science Program of the NRF of Korea (NRF-2022R1C1C1010052). J.-H. J. participated the introductory school of AGATES in Warsaw (Poland) and thanks the organizers for providing a good research environment throughout the school. The authors thank Hyun-Min Kim for invaluable advice and constant encouragement. The authors also thank Kangjin Han and Hayoung Choi for helpful discussion. This research was performed using the high-performance server computer provided by Finance-Fishery-Manufacture Industrial Mathematics Center on Big Data (FFMIMC). We would like to express our appreciation for this support.

\section{Preliminaries}\label{Sect:Preliminaries}
In this section, we review some conventions, definitions, and well-known facts we will frequently use in this paper. 

\subsection{Tensor Rank}
\begin{ntn} Throughout the paper, we use the following notations:
	\begin{itemize}
		\item $\mathbb{K}$ : a field of characteristic $\neq 2$;
		\item $V, V_i, W$ : finite dimensional $\mathbb{K}$-vector spaces;
		\item $V^{\ast}$ : the dual vector space of $V$;
		\item $[d] = \{1, 2, \cdots, d \}$ where $d$ is a positive integer.
	\end{itemize}
\end{ntn}

We follow the definitions and conventions in \cite{landsberg2012tensors, bernardi2018hitchhiker}.

\begin{dfn}[Multilinear map and tensor]\label{DefTensor}
	A map $\varphi:V_1 \times V_2 \times \cdots \times V_d \rightarrow W$ is said to be \emph{multilinear} if it is linear with respect to each $V_i$ for $i \in [d]$. The space of these multilinear maps is identified with $V_1^* \otimes V_2^* \otimes \cdots \otimes V_d^* \otimes W$. An element $\mathcal{T} \in V_1^* \otimes V_2^* \otimes \cdots \otimes V_d^* \otimes W$ is called a \emph{tensor}, and the number of vector spaces in the tensor product where $\mathcal T$ lives as an element is called the \emph{order of $\mathcal{T}$}.
\end{dfn}

 Since $V$ is of finite dimensional, the dual vector space $V^*$ is isomorphic to $V$. 
In case of $W=\mathbb{K}$ so that $W^{*} \simeq \mathbb{K}$, for instance $\operatorname{det}_n$ (on Section \ref{Subsec:Determinant}), tensoring $W$ or $W^*$ does not change the space $V_1^{*} \otimes \cdots \otimes V_d^{*}$. 
Hence, we regard $\mathcal{T} \in V_1 \otimes V_2 \otimes \cdots \otimes V_d \otimes W \cong V_1 \otimes V_2 \otimes \cdots \otimes V_d$ as an order $n$ tensor, not of order $n+1$. In addition, since a multilinear map $V_1 \times \cdots \times V_d \rightarrow W$ can be identified with a map $V_1 \times \cdots \times V_d \times W^* \rightarrow \mathbb{K}$, we deal with $V_1 \otimes V_2 \otimes \cdots \otimes V_d$ for further definitions.

\begin{dfn}[CP decomposition \cite{hitchcock1927expression, landsberg2012tensors, kolda2009tensor}]\label{DefCP_decomp}
	Let $\mathcal{T} \in V_1 \otimes V_2 \otimes \cdots \otimes V_d$ be a nonzero element. If there exist a positive integer $k$, for $i\in [k]$, and vectors $v_{i,j}\in V_{j}$ such that
	\begin{equation}\label{CP_decomp}
		\mathcal{T}=\sum_{i=1}^{k} v_{i,1} \otimes v_{i,2} \otimes \cdots \otimes v_{i,d}
	\end{equation}
	then \eqref{CP_decomp} is said to be a canonical polyadic decomposition or parallel factor decomposition of $\mathcal{T}$. Simply we call it CP decomposition of $\mathcal{T}$.
\end{dfn}

We can use this to define tensor rank as follows.

\begin{dfn}[Tensor rank]
	Let $\mathcal{T} \in V_1 \otimes V_2 \otimes \cdots \otimes V_d$ be a nonzero element. The number
	\begin{equation*}
		r=\min\left\{ k ~\middle|~\mathcal{T}=\sum_{i=1}^k v_{i,1} \otimes v_{i,2} \otimes \cdots \otimes v_{i,d}~~\text{where}~v_{i,j}\in V_{j}~~\text{for each}~j \in [d] \right\}
	\end{equation*}
	is called the \emph{(tensor) rank of $\mathcal{T}$}, and denoted by $\operatorname{rank}(\mathcal{T})$. 
\end{dfn}

It is well known that the rank is invariant under the change of bases and it generalizes the standard notion of the rank of matrices by considering a  matrix as an order $2$ tensor.
For example, the $2 \times 2$ determinant tensor  $\operatorname{det}_2$, which is defined as ${\det}_2=e_1^* \otimes e_2^* - e_2^* \otimes e_1^*$ for the standard dual basis $\{e_1^*,e_2^*\}$ of $\mathbb{K}^2$ (see Section \ref{Subsec:Determinant}), has a matrix representation
\begin{equation*}
	\begin{bmatrix}
		0 & -1 \\
		1 & 0
	\end{bmatrix}
\end{equation*}
which is of rank $2$, and thus $\operatorname{rank} ({\det}_2) = 2$. However, it is hard to determine the rank of a given tensor of order $d$ when $d>2$. Even finding an effective upper and lower bound of the tensor rank is a complicated question in most cases. In practice, the most direct method to improve an upper bound is to find a new decomposition formula using a lesser number of rank $1$ tensors. 

\subsection{SINDy}
SINDy is a popular data-driven method for identifying the governing equations of dynamical systems from time-series data. The method is based on the principle of sparsity, meaning that it seeks to identify the simplest set of equations that accurately describe the system's behavior\cite{kaiser2018sparse}. SINDy is also used in system modeling and control theory in various fields such as physics and chemistry\cite{brunton2016sparse}.

SINDy starts with the assumption that the dynamics of the system can be described by a set of nonlinear ordinary differential equations (ODEs) of the form:
\begin{equation}
	\frac{d\mathbf{x}}{dt} = \mathbf{f}(\mathbf{x})
\end{equation}
where $\mathbf{x} \in \mathbb{R}^n$ is the state vector of the system, $\mathbf{f}$ is a nonlinear function of the state, and $t$ is time. The goal of SINDy is to learn the shape of the function $\mathbf{f}$ directly from time-series data.

SINDy uses a sparse regression technique to identify the nonlinear terms in $\mathbf{f}$. To be precise, it tries to find the sparsest set of nonlinear terms that can be combined linearly to accurately describe the time evolution of the system. This can be formulated as an optimization problem
\begin{equation}\label{eqOptiPro}
	\underset{\boldsymbol{\theta}}{\operatorname{minimize}} \left( \left\|\mathbf{f}(\mathbf{x}) - \sum_{j=1}^p \theta_j \mathbf{\Phi}_j(\mathbf{x})\right\|_2^2 + \lambda \left\|\boldsymbol{\theta}\right\|_1 \right)
\end{equation}
where $\boldsymbol{\theta} \in \mathbb{R}^p$ is a vector of coefficients corresponding to the nonlinear terms, $\mathbf{\Phi}_j(\mathbf{x})$ is the $j$-th candidate function of the state, and $\lambda$ is a regularization parameter that controls the sparsity of the solution.

The first term in the (\ref{eqOptiPro}) represents the squared error between the true dynamics of the system and the linear combination of nonlinear terms. The second term is the same as the $l_1$-penalty term of LASSO, which improves sparsity of the solution by adding a penalty proportional to the absolute value of the coefficients. Such an optimization problem can be solved efficiently using numerical techniques such as gradient descent or proximal gradient methods. Once the coefficients $\boldsymbol{\theta}$ are found, the nonlinear terms can be combined linearly to obtain an estimate of the function $\mathbf{f}(\mathbf{x})$.

\subsection{SDL}
SDL is one of the popular techniques in the field of signal processing and machine learning. 
It is a technique that uses dictionary learning to effectively represent input data, where a dictionary is represented as a set of column vectors \cite{arora2015simple,rodriguez2018fast,kreutz2003dictionary, mairal2010online}.
SDL also uses the LASSO penalty to improve sparsity in the learned dictionary. 
We look for a set of basis vectors (or elements) that can efficiently represent a given set of signals (which is also used to update the dictionary, and thus we may expect a better representation the next time). 

SDL is defined as the following optimization problem
\begin{equation}\label{SDL}
	\underset{x}{\operatorname{minimize}} ~ \left( \frac{1}{2}\|y-D x\|_F^2+\lambda\|x\|_1 \right)
\end{equation}
where $\lambda$ is a regularization parameter. 
The above problem aims to keep the matrix $x$ sparse when representing a data matrix $y$ as a linear combination of column vectors (given by $x$) of a dictionary matrix $D$. Here, $\lambda$ is a regularization parameter that controls the sparsity of the solution as same as before.

SDL has several applications also in tensor analysis, for instance, Lee et al. propose a feature-sign search algorithm for finding a global optimum of the optimization problem \eqref{SDL} in a finite number of steps \cite{lee2006efficient}, and Duan et al. extended K-SVD \cite{aharon2006k} to the tensor scale and proposed the K-CPD algorithm for tensor sparse coding \cite{duan2012k}.

\section{Method}\label{Sect:Method}

SINDy and SDL both provide efficient methods for pursuing sparsity in data analysis and modeling.
Both methods require two steps, namely, building a set of candidates and then selecting an optimal solution that best represents the data.
We propose a numerical method to find a CP decomposition inspired by these two methods. 

For a given tensor $\mathcal{T}$, the process involves building a set of random rank $1$ tensors called the candidate set $D$. 
The aim is to find a sparse expression for the target $\mathcal T$ using as few candidates as possible.
If we regard $\mathcal T$ as the multilinear map, and if the codomain of $\mathcal T$ is one-dimensional, then the ordinary LASSO is sufficient, namely, we consider the following optimization problem

\begin{equation}\label{Method:LASSO}
	\textbf{x}_{sparse} =  \operatorname*{argmin}\limits_{X} \left( \|\textbf{y} - D\textbf{x}\|_{2} + \lambda\|\textbf{x} \|_1 \right) .
\end{equation}

We find a solution for $\mathbf{x}_{sparse}$ by minimizing the $l_2$-norm of the difference between the target tensor $\mathbf{y}$ and the product $D \mathbf{x}_{sparse}$, subject to a sparsity constraint $\lambda \|\textbf{x}\|_1$. 
The sparsity constraint is imposed using the $l_1$-norm regularization term. Note that the parameter $\lambda$ controls the trade-off between fitting the data and sparsity. The process is summarized in the following Figure \ref{Figure:LASSO1}.
\begin{figure}[h]
	\includegraphics[width=0.6\textwidth]{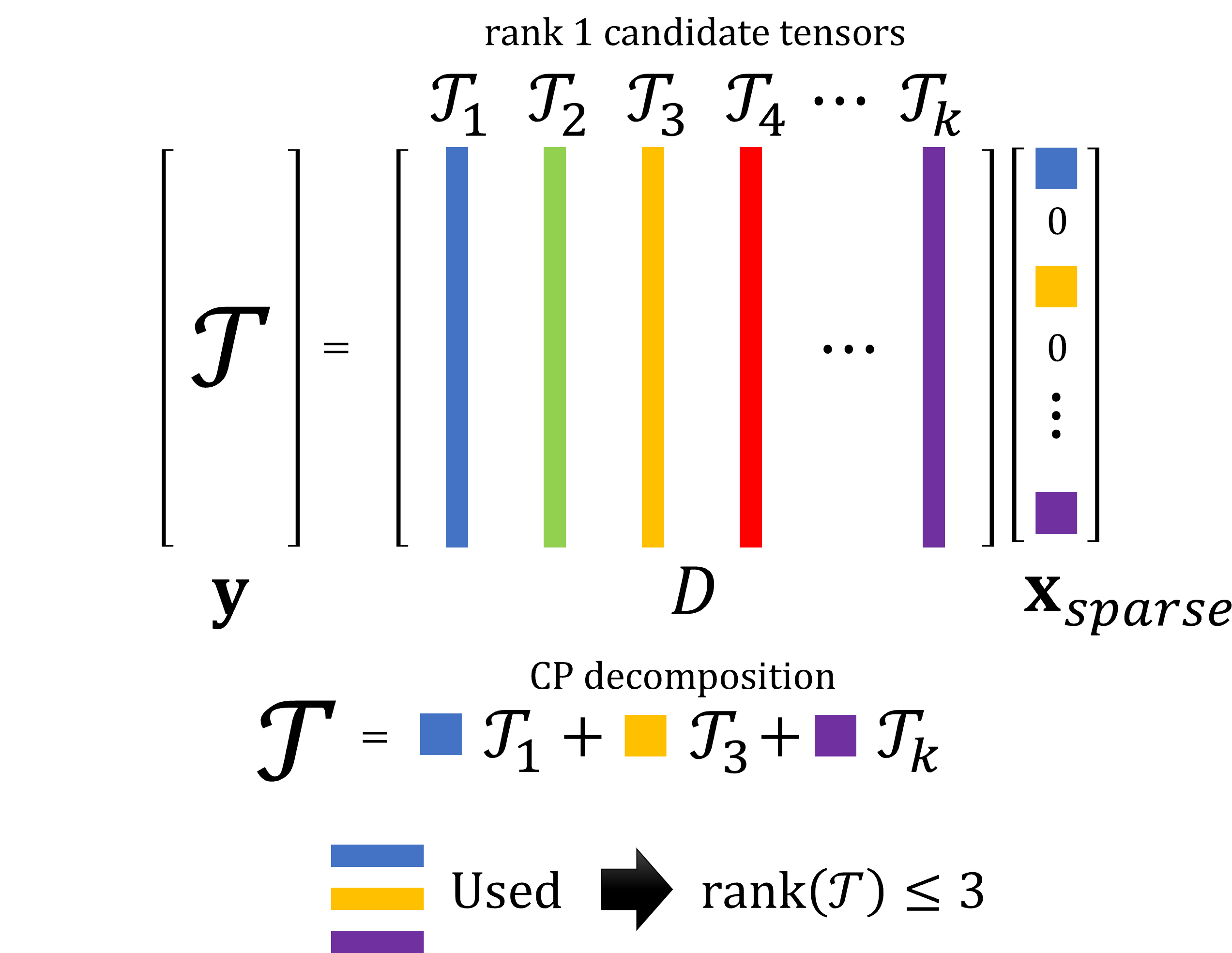}
	\caption{
		Schematic of calculating CP decompositions when the codomain of a tensor is $1$-dimensional. 
		Given $k$ candidate rank $1$ tensors, we build $\mathbf{y}$ from $n>k$ random data evaluated in the tensor $\mathcal{T}$, and an $n \times k$ matrix $D$ evaluated in each of the $k$ rank $1$ tensors $\mathcal T_i$.
		Solve equation \eqref{Method:LASSO} to find $\mathbf{x}_{sparse}$ that reconstructs $\mathbf{y}=D\mathbf{x}$.
		In this figure, $\mathcal{T}$ is decomposed by using three rank $1$ candidate tensors $\mathcal{T}_{1},\mathcal{T}_{3}$, and $\mathcal{T}_{k}$. 
		Hence, the upper bound of $\operatorname{rank}(\mathcal{T})$ is $3$.}\label{Figure:LASSO1}
\end{figure}

\begin{figure}[h]
	\includegraphics[width=0.65\textwidth]{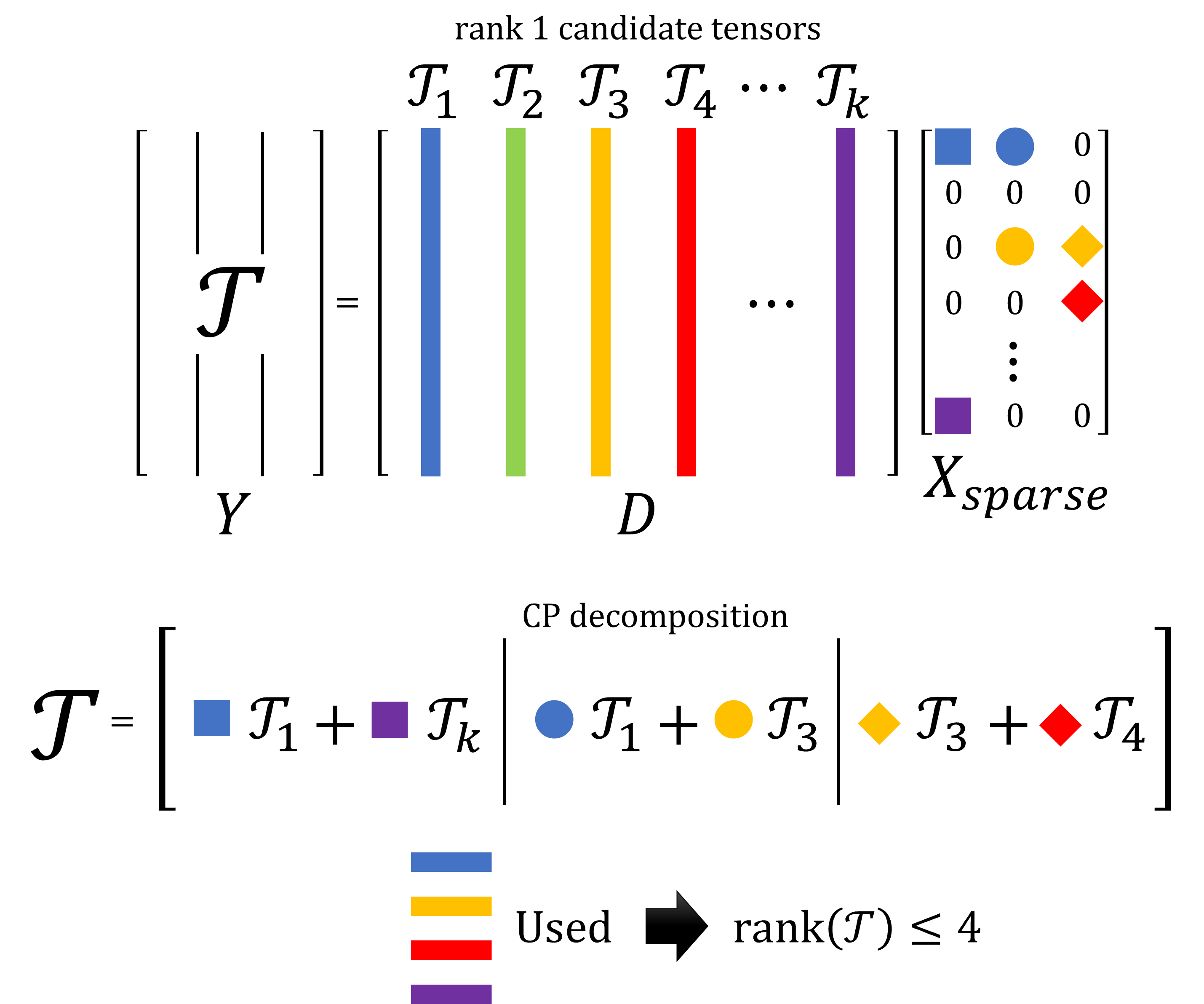}
	\caption{
		Schematic of calculating CP decompositions when the codomain of a tensor is $d(\geq 2)$-dimensional. Given $k$ candidate rank $1$ tensors, we build $Y$ from $n>k$ random data evaluated in the tensor $\mathcal{T}$, and an $n \times k$ matrix $D$ evaluated in each of the $k$ rank $1$ tensors $\mathcal T_i$.
		Solve equation \eqref{Method:LASSO2} to find the $k\times d$ matrix $X$ that reconstructs $Y=DX$.
		In this figure, $\mathcal{T}$ (maps to a $3$-dimensional space) is decomposed using six rank $1$ candidate tensors, but $\mathcal{F}(X_{sparse}) = 4$. We used four rank $1$ tensors $\mathcal{T}_{1},\mathcal{T}_{3},\mathcal{T}_{4}$, and $\mathcal{T}_{k}$. 
		Therefore, an upper bound of $\operatorname{rank}(\mathcal{T})$ is $4$.} \label{Figure:LASSO2}
\end{figure}

When the codomain of the multilinear map $\mathcal T$ has dimension at least $2$, we want to find a CP decomposition that minimizes the total number of candidate rank $1$ tensors used to represent the tensor.
This can be expressed as a group sparse optimization problem:
\begin{equation}\label{Method:LASSO2}
	X_{sparse} = \operatorname*{argmin}\limits_{X} \left(\|Y - DX\|_{F} + \lambda \mathcal{F} (X)\right). 
\end{equation}
Several penalty terms have been proposed to solve group sparse optimization problems, such as mixed norm\cite{wainwright2014structured,jiang2018joint}. 
However, since we are only interested in the number of rank $1$ candidate tensors used in the whole context, we may propose the following penalty term $\mathcal{F} (\cdot)$ for instance:
\begin{equation}\label{Group_penalty}
	\mathcal{F} (X) =\sum_{i=1}^{m}\delta\left(\sum_{j=1}^{n}\left|X_{i j}\right|\right)\quad \text{where}\quad \delta(x) = \begin{cases}
		1 & x\neq0\\
		0 & x=0
	\end{cases}
\end{equation}
If an inner summation $\sum_{j=1}^n |X_{ij}|$ is nonzero, then the corresponding row contains at least one nonzero element, and then the outer summation counts it towards the total number of nonzero rows. In other words, $\mathcal{F}(X)$ represents the number of nonzero rows in $X$.
Once the sparse solution $X_{sparse}$ is obtained, the number of nonzero rows in $X_{sparse}$ gives the number of candidates used in the sparse expression of the target tensor, see Figure \ref{Figure:LASSO2} for the summarization. As a result, the above optimization problem will track a more efficient CP decomposition of the target tensor as well as in the $1$-dimensional case. 

Big advantages of our method are that we can get a better CP decomposition if $D$ contains a sufficiently large class of rank $1$ candidate tensors, and then we can find a candidate CP decomposition of the given tensor efficiently via numerical approaches. 

\section{Experimental Result}\label{Sect:Result}
We provide some experiments and examples how to find effective CP decompositions for matrix multiplication and determinant tensors. 
\subsection{Matrix multiplications}
Let us explain our method with the example of matrix multiplication tensor $M_{\langle 2 \rangle}$. We follow the definition of the matrix multiplication tensor in \cite{landsberg2017geometry, fawzi2022discovering}.
\begin{dfn}
	The matrix multiplication tensor $M_{\langle n,m,p \rangle}$ is the multilinear map from $\mathbb{K}^{n \times m} \times \mathbb{K}^{m \times p} \rightarrow \mathbb{K}^{n \times p}$ defined by $(A,B) \mapsto AB$ (matrix multiplication). We denote $M_{\langle n \rangle}$ instead of $M_{\langle n,n,n \rangle}$.
\end{dfn}

What we consider is $M_{\langle 2 \rangle}$. Let $e_{ij}$ denote the $2 \times 2$ matrix of which the $(i,j)$-entry is $1$ and the others are $0$, so that $\{e_{11},e_{12},e_{21},e_{22}\}$ is a basis of $\mathbb{K}^{2 \times 2}$. Strassen's algorithm \cite{strassen1969gaussian} implies that $M_{\langle 2 \rangle}$ is represented as
\begin{equation}\label{Strassen}
	\begin{aligned}
		M_{\langle 2\rangle}=&(e_{11}^*+e_{22}^*)\otimes (e_{11}^* + e_{22}^*)\otimes (e_{11}+e_{22})\\
		&~+(e_{21}^*+e_{22}^*)\otimes e_{11}^* \otimes (e_{21}-e_{22})\\
		&~+e_{11}^*\otimes(e_{12}^*-e_{22}^*)\otimes(e_{12}+e_{22})\\
		&~+e_{22}^*\otimes(-e_{11}^*+e_{21}^*)\otimes (e_{11}+e_{21})\\
		&~+(e_{11}^*+e_{12}^*)\otimes e_{22}^* \otimes (-e_{11}+e_{12})\\
		&~+(-e_{11}^*+e_{21}^*)\otimes (e_{11}^*+e_{12}^*)\otimes e_{22}\\
		&~+(e_{12}^*-e_{22}^*)\otimes (e_{21}^*+e_{22}^*)\otimes e_{11},
	\end{aligned}
\end{equation}

where $\{e_{11}^*,e_{12}^*,e_{21}^*,e_{22}^*\}$ is a standard dual basis of $\mathbb{K}^{2 \times 2}$.

We may consider the matrix multiplication tensor $M_{\langle 2 \rangle}$ as a $2 \times 2$ matrix of homogeneous polynomials of degree $2$ with $8$ independent variables
\begin{equation*}
	\begin{bmatrix}
		A_{1,1} & A_{1,2}\\
		A_{2,1} & A_{2,2}
	\end{bmatrix}\begin{bmatrix}
		B_{1,1} & B_{1,2}\\
		B_{2,1} & B_{2,2}
	\end{bmatrix}=\begin{bmatrix}
		A_{1,1}B_{1,1}+A_{1,2}B_{2,1} & A_{1,1}B_{1,2}+A_{1,2}B_{2,2}\\
		A_{2,1}B_{1,1}+A_{2,2}B_{2,1} & A_{2,1}B_{1,2}+A_{2,2}B_{2,2}
	\end{bmatrix}.
\end{equation*}
We are going to compare CP decompositions for different candidate rank 1 tensors.

For $N>16$ random matrices $A$, $B$ and $C=AB$, we construct the candidate $\rank 1$ tensor $D$ and target tensor $Y$ as follows
\begin{align*}
	& D=
	\begin{bmatrix}
		| & | &   & | & |\\
		A_{1,1}B_{1,1} & A_{1,1}B_{1,2} & \cdots & A_{2,2}B_{2,1} & A_{2,2}B_{2,2}\\
		| & | &   & | & |
	\end{bmatrix} \in \mathbb{R}^{N \times 16}, \\
	& Y=
	\begin{bmatrix}
		| & | & | & |\\
		C_{1,1} & C_{2,1} & C_{1,2} & C_{2,2}\\
		| & | & | & |
	\end{bmatrix} \in \mathbb{R}^{N\times 4}.
\end{align*}

Since the tensor we want to decompose maps to a vector space of dimension $4$, we find a sparse solution $X_{sparse}$ through \eqref{Method:LASSO2}. From a sparse solution, we can read off a decomposition of the matrix multiplication tensor as the summand of the rank $1$ candidate tensor as follows.
\begin{align*}
	\begin{array}{cc}
		\begin{array}{c}
			p_1 = A_{1,1}B_{1,1},\\
			p_2 = A_{1,1}B_{1,2},\\
			p_3 = A_{1,2}B_{2,1},\\
			p_4 = A_{1,2}B_{2,2},\\
			p_5 = A_{2,1}B_{1,1},\\
			p_6 = A_{2,1}B_{1,2},\\
			p_7 = A_{2,2}B_{2,1},\\
			p_8 = A_{2,2}B_{2,2},
		\end{array}
		&
		\begin{array}{c}
			C_{1,1} = p_1 + p_3,\\
			C_{2,1} = p_5 + p_7,\\
			C_{1,2} = p_2 + p_4,\\
			C_{2,2} = p_6 + p_8.
		\end{array}
	\end{array}
\end{align*}

Note that this is the standard formula for the classical multiplication of between $2\times 2$ matrices. 
If we take a another $D$ so that it contains different rank $1$ candidate tensors, we expect a new decomposition of $M_{\langle 2 \rangle}$.
For $n$ random matrices $A,B,C=AB$, let vectors $A', B'$ be defined as follows
\begin{align*}&
	A' =
	\begin{bmatrix}
		A_{1,1}\\
		A_{2,2}\\
		A_{1,1} + A_{2,2}\\
		A_{2,1} + A_{2,2}\\
		A_{1,1} + A_{1,2}\\
		A_{1,1} - A_{2,1}\\
		A_{1,2} - A_{2,2}
	\end{bmatrix}, \quad
	B' =
	\begin{bmatrix}
		B_{1,1}\\
		B_{2,2}\\
		B_{1,1} + B_{2,2}\\
		B_{2,1} + B_{2,2}\\
		B_{1,1} + B_{1,2}\\
		B_{1,1} - B_{2,1}\\
		B_{1,2} - B_{2,2}
	\end{bmatrix}\in \mathbb{R}^{7}, 	
\end{align*}
so we allow sums and differences of some entries of $A$ and $B$. We construct rows of $D$ as the product of the elements of $A'$ and $B'$, respectively, and the rows of $Y$ by flattening the elements of $C$ by a choice of $N>7^2$ random matrices. 
Hence, each row of the matrices $D$ and $Y$ is of the form
\begin{align*}
	& D = 
	\begin{bmatrix}
		| & | & & | & |\\
		A'_{1} B'_{1} & A'_{1} B'_{2} & \cdots & A'_{7} B'_{6} & A'_{7} B'_{7}\\
		| & | & & | & |\\
	\end{bmatrix} \in \mathbb{R}^{N \times 7^2}\\
	& Y =
	\begin{bmatrix}
		| & | & | & |\\
		C_{1,1} & C_{2,1} & C_{1,2} & C_{2,2}\\
		| & | & | & |
	\end{bmatrix} \in \mathbb{R}^{N\times 4}
\end{align*}
for randomly chosen matrices $A, B$. For instance, $A'_{6} B'_{5}$ means the value corresponding to $(A_{1,1} - A_{2,1})(B_{1,1} + B_{1,2})$.
As similar as above, we find a sparse solution $X_{sparse}$ through \eqref{Method:LASSO2} with these new matrices. The nonzero rows of $X_{sparse}$ and the corresponding rank $1$ candidate tensor are given as follows.
\begin{equation*}
	\begin{array}{cc}
		\begin{array}{l}
			m_1=A_{1,1}(B_{1,2} - B_{2,2}),\\
			m_2=A_{2,2}(B_{1,1} - B_{2,1}),\\
			m_3=(A_{1,1} + A_{2,2})(B_{1,1} + B_{2,2}),\\
			m_4=(A_{2,1} + A_{2,2})B_{1,1},\\
			m_5=(A_{1,1} + A_{1,2})B_{2,2},\\
			m_6=(A_{1,1} - A_{2,1})(B_{1,1} + B_{1,2}),\\
			m_7=(A_{1,2} - A_{2,2})(B_{2,1} + B_{2,2}),
		\end{array}&
		\begin{array}{l}
			C_{1,1} = -m_2 + m_3 - m_5 + m_7,\\
			C_{2,1} = -m_2 + m_4,\\
			C_{1,2} = m_1 + m_5,\\
			C_{2,2} = m_1 + m_3 -m_4 - m_6.
		\end{array}
	\end{array}
\end{equation*}

Here we used only seven $\rank 1$ tensors $m_1, \ldots, m_7$ to represent the $2\times 2$ matrix multiplication, and it is easy to see that this is equivalent to the Stassen's algorithm \eqref{Strassen}.

\subsection{Determinants}\label{Subsec:Determinant}

Let us apply our method to find some new decompositions of determinant tensors.  We follow the definition of the determinant tensor in \cite{derksen2016nuclear,krishna2018tensor,houston2023new}.
\begin{dfn}
	Let $V$ be a vector space of dimension $n$. The Cartesian product $\underbrace{V \times \cdots \times V}_{n\text{-copies}}$ can be regarded as the space of $n \times n$ square matrices over $\mathbb{K}$. The $n \times n$ determinant tensor is the multilinear map which is defined as
	\begin{equation}\label{detn}
		{\det}_n= \sum_{\sigma \in S_n} sgn(\sigma) ~ e_{\sigma(1)}^* \otimes e_{\sigma(2)} ^*\otimes \cdots \otimes e_{\sigma(n)}^*
	\end{equation}
	where $\{e_1^*,...,e_n^*\}$ is a basis of $V^*$ and $sgn(\sigma)$ denotes the sign of the permutation $\sigma$. In particular,  $\det_n$ is an $n$-linear function from $\underbrace{V \times \cdots \times V}_{n\text{-copies}}$ to $\mathbb{K}$.
\end{dfn}

Note that the determinant $\det_n$ also can be regarded as a homogeneous polynomial of degree $n$ with $n^2$ independent variables parametrizing the entries of a square matrix. For example, ${\det}_2$ can be understood as a quadratic polynomial $x_{1,1}x_{2,2}-x_{1,2}x_{2,1}$ by considering the matrix $\begin{bmatrix} x_{1,1} & x_{1,2} \\ x_{2,1} & x_{2,2} \end{bmatrix}$ of independent variables. In this manner, a CP-decomposition of $\det_n$ as tensors provides an expression of $\det_n$ as homogeneous polynomials in terms of the sum of products of linear polynomials. 

\subsubsection{$3\times 3$ cases}

Let us recall some known results about the tensor rank of the determinant tensor ${\det}_n$. 
In \cite{derksen2016nuclear}, Derksen brought an explicit formula for the $3 \times 3$ matrix $A$, the determinant of $A$ is following:
\begin{equation}\label{Derk}
	\begin{aligned}
		{\det}_3(A) =\frac{1}{2}\big(& (A_{3,1}+A_{3,2}) (A_{2,1}-A_{2,2})  (A_{1,3}+A_{2,3})\\
		&+(A_{1,1}+A_{2,1}) (A_{2,2}-A_{3,2})  (A_{2,3}+A_{3,3})\\
		&+2A_{2,1}  (A_{3,2}-A_{1,2})  (A_{3,3}+A_{1,3}) \\
		&+(A_{3,1}-A_{2,1}) (A_{2,2}+A_{1,2})  (A_{2,3}-A_{1,3}) \\
		&+(A_{1,1}-A_{2,1})  (A_{3,2}+A_{2,2})  (A_{3,3}-A_{2,3})\big).
	\end{aligned}
\end{equation}
This formula implies that $\operatorname{rank}({\det}_3) \leq 5$ when $\operatorname{char}(\mathbb{K}) \neq 2$. 
Using (\ref{Derk}) and the generalized Laplace expansion, Derksen also showed that
\begin{equation}\label{UpperBound1} 
	\operatorname{rank}({\det}_n)\leq \left(\frac{5}{6} \right)^{\lfloor \frac{n}{3} \rfloor} n!,~~~~\text{if}~~\operatorname{char}(\mathbb{K}) \neq 2.
\end{equation}

Note that each term in the formula (\ref{Derk}) is a product of linear polynomials corresponding to an entry of the matrix $A$, or their sums and differences. Motivated by this observation, we consider the following matrix $M$
\begin{align*}
	M =
	\begin{bmatrix}
		A_{1,1} & A_{1,2} & A_{1,3}\\
		A_{2,1} & A_{2,2} & A_{2,3}\\
		A_{3,1} & A_{3,2} & A_{3,3}\\
		A_{1,1} - A_{2,1} & A_{1,2} - A_{2,2} & A_{1,3} - A_{2,3}\\
		A_{1,1} - A_{3,1} & A_{1,2} - A_{3,2} & A_{1,3} - A_{3,3}\\
		A_{2,1} - A_{3,1} & A_{2,2} - A_{3,2} & A_{2,3} - A_{3,3}\\
		A_{1,1} + A_{2,1} & A_{1,2} + A_{2,2} & A_{1,3} + A_{2,3}\\
		A_{1,1} + A_{3,1} & A_{1,2} + A_{3,2} & A_{1,3} + A_{3,3}\\
		A_{2,1} + A_{3,1} & A_{2,2} + A_{3,2} & A_{2,3} + A_{3,3}
	\end{bmatrix}
	\in \mathbb{R}^{9\times 3}
\end{align*}
for a randomly chosen matrix $A \in \mathbb{R}^{3 \times 3}$. 

We choose $N>9^3$ matrices $A\in \mathbb{R}^{3\times 3}$ randomly, and build the matrix
\begin{align*}
		& D =
		\begin{bmatrix}
			| & | &  & | & | \\
			M_{1,1}M_{1,2}M_{1,3} & M_{1,1}M_{1,2}M_{2,3} & \cdots & M_{3,1}M_{3,2}M_{2,3} & M_{3,1}M_{3,2}M_{3,3}\\
			| & | &  & | & |
		\end{bmatrix} \in \mathbb{R}^{N \times 9^3}.
\end{align*}

We consider the following linear equation
\begin{equation}\label{SINDy0}
	\mathbf{y}=\begin{bmatrix}
		\vline\\
		\det(A)\\
		\vline
	\end{bmatrix} = D\mathbf{x}
\end{equation}
to find an expression for $\det (A)$ as a linear combination of entries in $D$. Since the tensor we want to decompose maps to $1$-dimensional space, we obtain the sparse solution $\mathbf{x}_{sparse}$ of equation \eqref{SINDy0} through \eqref{Method:LASSO}.
The entries of $\mathbf{x}_{sparse}$ are the coefficients for each candidate rank $1$ tensor, and thus we read off them and obtain the formula for ${\det}_3$ as follows
\begin{equation}\label{three_det}
	\begin{aligned}
		{\det}_3 (A) = \frac{1}{2} 
		\big(&-2A_{1,1}(A_{2,2}+A_{3,2})(A_{2,2}-A_{3,2})\\
		&+(A_{1,1}-A_{2,1})(A_{1,2}-A_{3,2})(A_{1,2}+A_{3,2})\\
		&-(A_{1,1}-A_{3,1})(A_{1,2}-A_{2,2})(A_{1,2}+A_{2,2})\\
		&+(A_{1,1}+A_{2,1})(A_{1,2}+A_{3,2})(A_{1,2}-A_{3,2})\\
		&-(A_{1,1}+A_{3,1})(A_{1,2}+A_{2,2})(A_{1,2}-A_{2,2})\big).
	\end{aligned}
\end{equation}

It looks slightly different from Derksen's formula \eqref{Derk}, however, it is easy to see that the equation \eqref{three_det} also computes the determinant of a $3\times 3$ matrix with five similar terms.

\subsubsection{$4\times 4$ case}

We are going to find a CP decomposition of $\det_4$ using $16$ independent variables in the matrix of size $4\times 4$
\begin{equation*}
	A=\begin{bmatrix}
		A_{1,1} & A_{1,2} & A_{1,3} & A_{1,4}\\
		A_{2,1} & A_{2,2} & A_{2,3} & A_{2,4}\\
		A_{3,1} & A_{3,2} & A_{3,3} & A_{3,4}\\
		A_{4,1} & A_{4,2} & A_{4,3} & A_{4,4}\\
	\end{bmatrix}.
\end{equation*}
We pay attention to terms in Derksen's formula \eqref{Derk} for $\det_3$ which are composed of products of the sums or differences of two entries in each row (column) of the matrix, or entries of the matrix. 

We construct the following matrix $M$ for $4\times 4$ matrix $A$ composed of the entries, sums, and differences of the entries of $A$.
Each element in the column of $M$ is filled with the value of an entry of $A$, or adding/subtracting another entry in the column of $A$. 
$$
M=
\begin{bmatrix}
	A_{1,1} & A_{1,2} & A_{1,3} & A_{1,4}\\
	A_{2,1} & A_{2,2} & A_{2,3} & A_{2,4}\\
	A_{3,1} & A_{3,2} & A_{3,3} & A_{3,4}\\
	A_{4,1} & A_{4,2} & A_{4,3} & A_{4,4}\\ 
	A_{1,1}-A_{2,1} & A_{1,2}-A_{2,2} & A_{1,3}-A_{2,3} & A_{1,4}-A_{2,4}\\
	A_{1,1}-A_{3,1} & A_{1,2}-A_{3,2} & A_{1,3}-A_{3,3} & A_{1,4}-A_{3,4}\\
	A_{1,1}-A_{4,1} & A_{1,2}-A_{4,2} & A_{1,3}-A_{4,3} & A_{1,4}-A_{4,4}\\
	A_{2,1}-A_{3,1} & A_{2,2}-A_{3,2} & A_{2,3}-A_{3,3} & A_{2,4}-A_{3,4}\\
	A_{2,1}-A_{4,1} & A_{2,2}-A_{4,2} & A_{2,3}-A_{4,3} & A_{2,4}-A_{4,4}\\
	A_{3,1}-A_{4,1} & A_{3,2}-A_{4,2} & A_{3,3}-A_{4,3} & A_{3,4}-A_{4,4}\\
	A_{1,1}+A_{2,1} & A_{1,2}+A_{2,2} & A_{1,3}+A_{2,3} & A_{1,4}+A_{2,4}\\
	A_{1,1}+A_{3,1} & A_{1,2}+A_{3,2} & A_{1,3}+A_{3,3} & A_{1,4}+A_{3,4}\\
	A_{1,1}+A_{4,1} & A_{1,2}+A_{4,2} & A_{1,3}+A_{4,3} & A_{1,4}+A_{4,4}\\
	A_{2,1}+A_{3,1} & A_{2,2}+A_{3,2} & A_{2,3}+A_{3,3} & A_{2,4}+A_{3,4}\\
	A_{2,1}+A_{4,1} & A_{2,2}+A_{4,2} & A_{2,3}+A_{4,3} & A_{2,4}+A_{4,4}\\
	A_{3,1}+A_{4,1} & A_{3,2}+A_{4,2} & A_{3,3}+A_{4,3} & A_{3,4}+A_{4,4}
\end{bmatrix}
$$
 
We expect that $\det_4$ also can be expressed by a formula in a similar shape of terms performing Derksen's formula \eqref{Derk} for $\det_3$. We take one row of the candidates (which will be a row of the matrix $D$) as follows.
$$
\underbrace{
	\begin{bmatrix}
		M_{1,1}M_{1,2}M_{1,3}M_{1,4} & M_{1,1}M_{1,2}M_{1,3}M_{2,4} &\cdots& 
		M_{16,1}M_{16,2}M_{16,3}M_{16,4}
\end{bmatrix}}_{16^4\times 1 \text{ vector}}
$$
Here, the entries are the product of $4$ entries by picking an entry from each column of $M$, for example, $M_{6,1}M_{11,2}M_{1,3}M_{2,4}$ represents $(A_{1,1}-A_{3,1})(A_{1,2}+A_{2,2})(A_{1,3})(A_{2,4})$. Since $M$ contains products of the entries $A_{i,j}$ of $A$, it is clear that $\det_4$ can be written as a linear combination of entries in $M$. The length of the candidate vector is $16^{4}$, hence, to obtain an explicit relation among them, we construct the following matrix $D$ built from $4\times 4$ random matrices $A$ more than $16^4$. 
$$
D = 
\begin{bmatrix}
	\vline & \vline &  & \vline\\
	M_{1,1}M_{1,2}M_{1,3}M_{1,4} & M_{1,1}M_{1,2}M_{1,3}M_{2,4} &\cdots& 
	M_{16,1}M_{16,2}M_{16,3}M_{16,4}\\
	\vline & \vline &  & \vline\\
\end{bmatrix} \in \mathbb{R}^{N\times 16^4}
$$

We consider the following linear equation
\begin{equation}\label{SINDy1}
	\mathbf{y}=\begin{bmatrix}
		\vline\\
		\det(A)\\
		\vline
	\end{bmatrix} = D\mathbf{x}
\end{equation}
and find a sparse solution $\mathbf{x}$ of the sparse optimization problem through LASSO \eqref{Method:LASSO} which minimizes
$$
\mathbf{x}_{\text{sparse}} = \underset{\mathbf{x}}{\operatorname{argmin}} \big(\|\mathbf{y} - D\mathbf{x} \|_{2} +\lambda \|\mathbf{x} \|_{1}\big).
$$
A sparse solution $\mathbf{x}_{\text{sparse}}$ will give the coefficient of candidates in $D$, and thus describes $\det_4$ as a linear combination of the basic elements $M_{i,1} M_{i,2} M_{i,3} M_{i,4}$. 
Indeed, we obtain the following new determinant formula
\begin{equation*}
	\begin{aligned}
		{\det}_4(A) = \frac{1}{2}\big(
		&(A_{1,1}-A_{2,1})(A_{3,2}-A_{4,2})(A_{3,3}+A_{4,3})(A_{1,4}+A_{2,4})\\
		&-(A_{1,1}-A_{3,1})(A_{2,2}-A_{4,2})(A_{2,3}+A_{4,3})(A_{1,4}+A_{3,4})\\
		&+(A_{1,1}-A_{4,1})(A_{2,2}-A_{3,2})(A_{2,3}+A_{3,3})(A_{1,4}+A_{4,4})\\
		&+(A_{2,1}-A_{3,1})(A_{1,2}-A_{4,2})(A_{1,3}+A_{4,3})(A_{2,4}+A_{3,4})\\
		&-(A_{2,1}-A_{4,1})(A_{1,2}-A_{3,2})(A_{1,3}+A_{3,3})(A_{2,4}+A_{4,4})\\
		&+(A_{3,1}-A_{4,1})(A_{1,2}-A_{2,2})(A_{1,3}+A_{2,3})(A_{3,4}+A_{4,4})\\
		&+(A_{1,1}+A_{2,1})(A_{3,2}+A_{4,2})(A_{3,3}-A_{4,3})(A_{1,4}-A_{2,4})\\
		&-(A_{1,1}+A_{3,1})(A_{2,2}+A_{4,2})(A_{2,3}-A_{4,3})(A_{1,4}-A_{3,4})\\
		&+(A_{1,1}+A_{4,1})(A_{2,2}+A_{3,2})(A_{2,3}-A_{3,3})(A_{1,4}-A_{4,4})\\
		&+(A_{2,1}+A_{3,1})(A_{1,2}+A_{4,2})(A_{1,3}-A_{4,3})(A_{2,4}-A_{3,4})\\
		&-(A_{2,1}+A_{4,1})(A_{1,2}+A_{3,2})(A_{1,3}-A_{3,3})(A_{2,4}-A_{4,4})\\
		&+(A_{3,1}+A_{4,1})(A_{1,2}+A_{2,2})(A_{1,3}-A_{2,3})(A_{3,4}-A_{4,4})\big).
	\end{aligned}
\end{equation*}
We can check by direct comparison with the explicit $\det_4 (A)$ that the above formula is valid. From this formula, we also have the following CP decomposition for $\det_4$:
\begin{thm}\label{MainThm} 
Let $char(\mathbb{K}) \neq 2$. Then
\begin{equation}\label{newdet4}
\begin{aligned}
	{\det}_4=\frac{1}{2}\big( &(e^{*}_1-e^{*}_2)\otimes (e^{*}_3-e^{*}_4)\otimes (e^{*}_3+e^{*}_4)\otimes (e^{*}_1+e^{*}_2)\\
		&-(e^{*}_1-e^{*}_3)\otimes (e^{*}_2-e^{*}_4)\otimes (e^{*}_2+e^{*}_4)\otimes (e^{*}_1+e^{*}_3)\\
		&+(e^{*}_1-e^{*}_4)\otimes (e^{*}_2-e^{*}_3)\otimes (e^{*}_2+e^{*}_3)\otimes (e^{*}_1+e^{*}_4)\\
		&+(e^{*}_2-e^{*}_3)\otimes (e^{*}_1-e^{*}_4)\otimes (e^{*}_1+e^{*}_4)\otimes (e^{*}_2+e^{*}_3)\\
		&-(e^{*}_2-e^{*}_4)\otimes (e^{*}_1-e^{*}_3)\otimes (e^{*}_1+e^{*}_3)\otimes (e^{*}_2+e^{*}_4)\\
		&+(e^{*}_3-e^{*}_4)\otimes (e^{*}_1-e^{*}_2)\otimes (e^{*}_1+e^{*}_2)\otimes (e^{*}_3+e^{*}_4)\\
		&+(e^{*}_1+e^{*}_2)\otimes (e^{*}_3+e^{*}_4)\otimes (e^{*}_3-e^{*}_4)\otimes (e^{*}_1-e^{*}_2)\\
		&-(e^{*}_1+e^{*}_3)\otimes (e^{*}_2+e^{*}_4)\otimes (e^{*}_2-e^{*}_4)\otimes (e^{*}_1-e^{*}_3)\\
		&+(e^{*}_1+e^{*}_4)\otimes (e^{*}_2+e^{*}_3)\otimes (e^{*}_2-e^{*}_3)\otimes (e^{*}_1-e^{*}_4)\\
		&+(e^{*}_2+e^{*}_3)\otimes (e^{*}_1+e^{*}_4)\otimes (e^{*}_1-e^{*}_4)\otimes (e^{*}_2-e^{*}_3)\\
		&-(e^{*}_2+e^{*}_4)\otimes (e^{*}_1+e^{*}_3)\otimes (e^{*}_1-e^{*}_3)\otimes (e^{*}_2-e^{*}_4)\\
		&+(e^{*}_3+e^{*}_4)\otimes (e^{*}_1+e^{*}_2)\otimes (e^{*}_1-e^{*}_2)\otimes (e^{*}_3-e^{*}_4)\big).
\end{aligned}
\end{equation}
\end{thm}

In particular, we have $\operatorname{rank}({\det}_4) \leq 12$ when $\operatorname{char}(\mathbb{K}) \neq 2$. Very recently, Houston et al. found a formula of ${\det}_4$ which consists of $15$ summands in case of $\operatorname{char}(\mathbb{K}) \neq 2$, and a formula which consists of $12$ summands in case of $\operatorname{char}(\mathbb{K})=2$ \cite{houston2023new}. By combining these results and ours, we have $\operatorname{rank}({\det}_4) \leq 12$ over an arbitrary field.

In addition, by considering $\det_4$ as a homogeneous polynomial of degree $4$ with $16$ independent variables and denoting it ${\det}_4$ again, we also have $\operatorname{Wrank}(det_4) \leq 2^3 \cdot \operatorname{rank}(det_4) \leq 96$ when $\operatorname{char}(\mathbb{K}) \neq 2,3$ \cite{ranestad2011rank}. Here, $\operatorname{Wrank}$ implies the Waring rank of a symmetric tensor (or equivalently, a homogeneous polynomial) \cite{landsberg2012tensors}. This inequality achieves the same bound as $\operatorname{Wrank}({\det}_n)\leq n \cdot n!$  \cite{johns2022improved} for $n=4$ in a different way.

We finish by a remark that the formula \eqref{newdet4} encodes a lot of symmetries, in particular, on the choice of partitions of $[n]$ into subsets of size $2$. This gives rise to a general formula for the determinant of $n \times n$ matrices for arbitrary $n$, see \cite{jkk2023detn}.

\bibliography{Det}
\vspace{0.5cm}

\end{document}